\documentclass[11pt]{article}

\usepackage{amssymb,amsmath,euscript}
\makeatletter \@addtoreset{equation}{section} \makeatother
\newtheorem{proposition}{Proposition}[section]
\newtheorem{lemma}[proposition]{Lemma}
\newtheorem{theorem}[proposition]{Theorem}

\newtheorem{remark}[proposition]{Remark}

\def\R{{\mathbb R}}

\def\E{{\mathbb E}}
\def\F{{\mathcal F}}
\def\S{{\mathcal S}}

\def\P{{\mathbb P}}

\def\1{\mathbf{1}}
\def\a{\alpha}
\def\la{\lambda}

\def\de{\delta}
\def\ga{\gamma}

\def\th{\theta}
\def\si{\sigma}

\def\vf{\varphi_q}

\def\cqfd{$\square$}

\begin{document}
\title{Sharp estimates on the tail behavior of a multistable distribution}
\author{Antoine Ayache\\UMR CNRS 8524, Laboratoire Paul Painlev\'e, B\^at. M2\\
  Universit\'e Lille 1\\ 59655 Villeneuve d'Ascq Cedex, France\\
E-mail: \texttt{Antoine.Ayache@math.univ-lille1.fr}
}
\maketitle

\begin{abstract}
Multistable distributions are natural extensions of symmetric $\a$ stable distributions. They have been introduced quite recently by Falconer, L\'evy V\'ehel and 
their co-authors in \cite{FGL09,FL09,FLiu10}. Roughly speaking such a distribution is obtained by replacing the constant parameter $\a$ of a symmetric stable distribution by a (Lebesgue) mesurable function 
$\a (x)$ with values in $[a,2]$, where $a>0$ denotes a fixed arbitrarily small real number. 

Let $Y$ be an arbitrary symmetric $\a$ stable random variable of scale parameter $\si>0$, an important classical result concerning the heavy-tailed behavior of its distribution (see e.g. \cite{ST94})
is that there exists an explicit constant $C(\a)>0$, only depending on $\a\in (0,2)$, such that $\lim_{\la\rightarrow +\infty} \big (C(\a)\si^{\a}\la^{-\a}\big)^{-1}\P\big(|Y|>\la\big)=1$. In this article, we show that the latter result can be extended to the setting of multistable random variables,
when the function $\a (x)$ is with values in an arbitray compact interval $[a,b]$ contained in $(0,2)$.
\end{abstract}

\section{Introduction and statement of the main result}
\label{sec:intro}
Let us first briefly recall the definition of the $\a (x)$-multistable stochastic integral which was quite recently introduced in \cite{FLiu10}. To this end, we need to fix some definitions and notations.
We denote by $a>0$ a fixed arbitrarily small real number and we denote by $\a:\R\rightarrow [a,2]$ an arbitrary fixed Lebesgue mesurable function; $\F_\a$, the corresponding variable exponent Lebesgue space, is defined as,
$$
\F_\a=\Big\{f:\mbox{ $f$ is Lebesgue mesurable with $\int_{\R} \big|f(x)\big|^{\a(x)}\,dx <\infty$}\Big\}.
$$
The space $\F_\a$ is equipped with a quasinorm denoted by $\|\cdot\|_\a$; for every $f\in \F_{\a}^*=\F_{\a}\setminus\{0\}$, $\|f\|_\a$ is defined as the unique 
$\la_{0,f}\in (0,+\infty)$, such that, $\int_{\R} \big|\la_{0,f}^{-1}f(x)\big|^{\a(x)}\,dx=1$, moreover one sets $\|0\|_\a=0$. Let us recall that to say that
$\|\cdot\|_\a$ is a quasinorm means that $\|\cdot\|_\a$ satisfies the following 3 properties:
\begin{itemize}
\item for all $f\in\F_\a$, one has $f=0$ if and only if $\|f\|_\a=0$;
\item for all $f\in\F_\a$ and $\de\in\R$, one has $\|\de f\|_\a =|\de|\|f\|_\a$;
\item there is a constant $k>0$, such that for all $f,g\in\F_\a$, one has $\|f+g\|_\a \le k\big ( \|f\|_\a+\|g\|_\a \big)$ (weak triangle inequality).
\end{itemize}
The following theorem is an important result of \cite{FLiu10}, which allows to define on $\F_\a$ the multistable stochastic integral; it has been obtained thanks to Kolmogorov's extension Theorem.

\begin{theorem}
\label{th:falliu}
\cite{FLiu10} There exists a real-valued stochastic process indexed by the space $\F_\a$, denoted by $\big\{I(f):f\in\F_\a\big\}$, whose finite dimensional distributions are characterized 
by the following property: for all integer $d\ge 1$ and all $f_1,\ldots, f_d\in \F_\a$, $\Phi_{I(f_1),\ldots,I(f_d)}=\Phi_{f_1,\ldots,f_d}$ the characteristic function of the random vector $\big (I(f_1),\ldots, I(f_d)\big)$, satisfies for all $(\th_1,\ldots,\th_d)\in\R^d$,
\begin{equation}
\label{eq1:falliu}
\Phi_{f_1,\ldots,f_d}(\th_1,\ldots,\th_d)=\exp\Big\{-\int_{\R}\big|\sum_{l=1}^d \th_l f_l(x)\big|^{\a(x)}\,dx\Big\}.
\end{equation}
\end{theorem}
Recall that, generally speaking, the distribution of an arbitrary random vector $(X_1,\ldots,X_d)$ is completely determined by $\Phi_{X_1,\ldots,X_d}$ its 
characteristic function, which is defined for all $(\th_1,\ldots,\th_d)\in\R^d$, as,
\begin{equation}
\label{eq:vectchar}
\Phi_{X_1,\ldots,X_d}(\th_1,\ldots,\th_d)=\E\Big (\exp i\sum_{l=1}^d \th_l X_l\Big).
\end{equation}
For each $f\in\F_\a$, $I(f)$ is called an $\a(x)$-multistable random variable and its distribution is called an $\a(x)$-multistable distribution. Generally speaking,
in many applied and theoretical problems, it is important to have a sharp estimates on the tail behavior of a probability distribution. The following theorem, which is our main result,
provides such an estimation in the case of an $\a(x)$-multistable distribution.
\begin{theorem}
\label{th:main}
Assume that there is $b\in (a,2)$ such that for almost all $x\in\R$, $\a(x)\in [a,b]$.
Let $C$ be the continuous strictly positive function defined for all $\ga\in [a,b]$ as,
\begin{equation}
\label{eq3:defC}
C(\ga)= \frac{2}{\pi} \mbox{ if $\ga=1$, and } C(\ga)=\frac{1-\ga}{\Gamma (2-\ga)\cos\big(2^{-1}\pi\ga\big)} \mbox{ else,}
\end{equation}
where $\Gamma$ is the usual "Gamma" function. For each $f\in\F_{\a}$ and real number $\la > 0$, let us set,
\begin{equation}
\label{eq1:limdomterm}
T_f (\la)=\int_\R \big |\la^{-1}f(x)\big |^{\a (x)} C\big(\a(x)\big)\, dx.
\end{equation}
Then, one has,  
\begin{equation}
\label{eq1:main}
\lim_{\la\rightarrow +\infty} \left\{\sup_{f\in\S_\a}\Big|\frac{\P\big(|I(f)|>\la\big)}{T_f(\la)}-1\Big|\right\}=0,
\end{equation}
where $\S_\a=\big\{f\in\F_\a: \|f\|_\a=1\big\}$ denotes the unit sphere of $\F_\a$.
\end{theorem}

Before ending this introduction, let us make some remaks concerning Theorem~\ref{th:main}.

\noindent {\bf Remarks:}
\begin{itemize}
\item Theorem~\ref{th:main} is an extension to the setting of multistable random variables of Property~1.2.15, on page 16 in \cite{ST94}. Indeed, assuming that for almost
all $x\in\R$, $\a (x)=\a$ where $\a\in (0,2)$ is a constant, then $I(f)$ reduces to a usual symmetric $\a$ stable random variable of scale parameter $\sigma=\Big (\int_{\R} 
\big |f(x)\big|^\a\,dx\Big )^{1/\a}$ and $T_f (\la)$ reduces to $C(\a)\sigma^{\a}\la^{-\a}$; thus we recover the statement of Property 1.2.15, on page 16 in \cite{ST94}.
\item Theorem~\ref{th:main} shows that when $\|\a\|_{L^\infty (\R)}<2$, then the distribution of the multistable random variable $I(f)$ is heavy-tailed (see 
(\ref{eq2:limdomterm})). 
\end{itemize}

\section{Proof of the main result}
\label{sec:proof}
The main goal of this section is to prove Theorem~\ref{th:main}. To this end, we need to introduce some notations and to derive some preliminary result.
We denote by $f$ an arbitrary function of $\F_{\a}^*$ and by $I(f)$ the $\a(x)$-multistable random variable defined as the $\a(x)$-multistable stochastic integral of $f$.
The characteristic function of $I(f)$ is denoted by $\Phi_f$, recall that it is defined for all $\th\in\R$ as,
\begin{equation}
\label{eq:defcarX}
\Phi_f (\th)=\E\big(e^{i\th I(f)}\big).
\end{equation}
Observe that, in view of (\ref{eq1:falliu}), one has for all $\th\in\R$,
\begin{equation}
\label{eq:charX}
\Phi_f(\th)=\exp\Big\{-\int_\R \big|\th f(x)\big|^{\a(x)}\,dx\Big\}.
\end{equation}
As a consequence:

\begin{remark}
\label{rem:carLp}
$\Phi_f$ is an even function which belongs to the Lebesgue space $L^p (\R)$, for any arbitrary $p\in (0,+\infty)$, in particular it belongs 
to $L^1 (\R)$. Therefore, the distribution of $I(f)$ is absolutely continuous with respect to the Lebesgue measure on $\R$, moreover
$D_f$ the corresponding Randon-Nikodym derivative (i.e. the probability density function of $I(f)$) is 
given for all $x\in\R$, by
\begin{equation}
\label{eq:defg}
D_f(x)=(2\pi)^{-1}\widehat{\Phi}(x)=(2\pi)^{-1}\int_{\R} e^{-ix\th} \Phi_f(\th)\,d\th,
\end{equation}
which implies that $D_f$ is a continuous, even and bounded function. Moreover, for all $\th\in\R$,
\begin{equation}
\label{eq:four-car}
\Phi_f(\th)=\widehat{D_f}(\th)=\int_\R e^{-i\th x} D_f (x)\,dx.
\end{equation}
\end{remark}

Notice that throughout this paper the Fourier transform of an arbitrary function $h$ of $L^1 (\R)$, is defined, for all $x\in\R$, as $\widehat{h}(x)=\int_\R e^{-ix\th} 
h(\th)\,d\th$.

\noindent {\sc Proof of Remark~\ref{rem:carLp}:} In view of (\ref{eq:charX}), it is clear that $\Phi_f$ is an even function, moreover by using the fact that for almost
all $x\in\R$, $\a(x)\in [a,b]\in (0,2)$, one has for all $p\in (0,+\infty)$,
\[
\begin{split}
& \int_\R \big|\Phi_f(\th)\big|^p\,d\th=\int_{\R} \exp\Big\{-p\int_\R \big|\th f(x)\big|^{\a(x)}\,dx\Big\}\,d\th\\
& \le \int_{|\th|\le 1} \exp\Big\{-p|\th|^b\int_\R \big|f(x)\big|^{\a(x)}\,dx\Big\}\,d\th+\int_{|\th|>1} \exp\Big\{-p|\th|^a\int_\R \big|f(x)\big|^{\a(x)}\,dx\Big\}\,d\th<\infty.
\end{split}
\]
\cqfd

Let $q$ be an arbitrary fixed real number strictly larger than $1$, and let $\vf$ be an even real-valued $C^\infty$ function whose Fourier transform $\widehat{\vf}$ is an even compactly supported $C^5$ function with values in $[0,1]$ satisfying for all $x\in\R$,
\begin{equation}
\label{eq:phihat}
\widehat{\vf}(x)=
\left\{
\begin{array}{l}
1 \mbox{ if $|x|\le 1$,}\\
0 \mbox{ if $|x| \ge \frac{1+q}{2}$.}
\end{array}
\right.
\end{equation}
Observe that for all $\gamma\in [0,4)$,
\begin{equation}
\label{eq:sobophi}
\int_{\R}\big (1+|\th|\big )^\ga \vf (\th)\,d\th<\infty.
\end{equation}
Also, observe that one has for all integer $j\ge 0$,
\begin{equation}
\label{eq:supphihat}
\widehat{\vf}(q^{-j}x)=0 \mbox{ when $|x|\ge \Big (\frac{1+q}{2}\Big ) q^j$,}
\end{equation}
and 
\begin{equation}
\label{eq:comsupphihat}
1-\widehat{\vf}(q^{-j}x)=0 \mbox{ when $|x|\le q^{j}$.}
\end{equation}
For all $\la\in [q,+\infty)$, let $j_0(\la,q)\ge 1$ be the unique integer such that
\begin{equation}
\label{eq:dej0}
q^{j_0(\la,q)}\le \la <q^{j_0(\la,q)+1},
\end{equation}
therefore, denoting by $[\cdot]$ the integer part function, it follows that,
\begin{equation}
\label{eq:expj0}
j_0(\la,q)=\Big [\frac{\log \la}{\log q}\Big ].
\end{equation}
\begin{lemma}
\label{lem:encad-prob}
For all real number $\xi\ge 1 $, let us set,
\begin{equation}
\label{eq1:minoretaj}
\eta_f (\xi)=\int_{\R}\vf (\th)\left(1-\exp\Big\{-\int_\R \big|\xi^{-1}\th f(x)\big|^{\a(x)}\,dx\Big\}\right)\,d\th.
\end{equation}
Then, for each $\la\in [q,+\infty)$ one has,
\begin{equation}
\label{eq:LPProbB}
\eta_f\big(q^{j_0(\la,q)+1}\big)\le \P\big(|I(f)|>\la\big)\le \eta_f \big(q^{j_0(\la,q)-1}\big).
\end{equation}
\end{lemma}

\noindent {\sc Proof of Lemma~\ref{lem:encad-prob}:} 
Using the fact that $D_f$ is the probability density function of $I(f)$, one has that 
$$
\P\big(|I(f)|>\la\big)=\int_{|x|>\la}D_f (x)\,dx.
$$ 
Therefore, it follows from (\ref{eq:dej0}) and (\ref{eq:comsupphihat}) that,
\begin{eqnarray}
\label{eq1:encad-prob}
\nonumber
\P\big(|I(f)|>\la\big)\ge \int_{|x|\ge q^{j_0(\la,q)+1}}D_f (x)\,dx &\ge & \int_{|x|\ge q^{j_0(\la,q)+1}}\big (1-\widehat{\vf}(q^{-j_0(\la,q)-1}x)\big) D_f (x)\,dx\\
&=& \int_{\R}\big (1-\widehat{\vf}(q^{-j_0(\la)-1}x)\big) D_f (x)\,dx;
\end{eqnarray}
on the other hand, (\ref{eq:dej0}) and (\ref{eq:supphihat}) imply that
\begin{eqnarray}
\label{eq2:encad-prob}
\nonumber
\P\big(|I(f)|>\la\big)&\le& \int_{|x|\ge q^{j_0(\la,q)}}D_f(x)\,dx\\
\nonumber
&=& \int_{|x|\ge q^{j_0(\la)}}\widehat{\vf}(q^{-j_0(\la,q)+1}x)D_f(x)\,dx+\int_{|x|\ge q^{j_0(\la)}}
\big (1-\widehat{\vf}(q^{-j_0(\la,q)+1}x)\big) D_f(x)\,dx\\
&\le & \int_{\R}\big(1-\widehat{\vf}(q^{-j_0(\la,q)+1}x)\big) D_f (x)\,dx.
\end{eqnarray}
Let us now prove that for all real number $\de >0$, one has,
\begin{equation}
\label{eq3:encad-prob}
\int_{\R}\big(1-\widehat{\vf}(\de x)\big) D_f(x)\,dx=\int_{\R}\vf (\th)\big (1-\widehat{D_f}(\de\th)\big)\,d\th.
\end{equation}
In view of the fact that 
$$
\int_\R D_f (x)\,dx=1 \mbox{ and } \int_\R \vf (\th)\,d\th=\widehat{\vf}(0)=1,
$$
it is sufficient to show that
$$
\int_{\R}\widehat{\vf}(\de x)D_f(x)\,dx=\int_{\R}\vf (\th)\widehat{D_f}(\de\th)\,d\th.
$$
By using the definition of the Fourier transform of an $L^1(\R)$ function and Fubini Theorem, it follows that
\begin{eqnarray*}
\int_{\R}\widehat{\vf}(\de x) D_f (x)\,dx&=& \int_{\R}\Big (\int_\R e^{-i\de x\th}\vf(\th)\,d\th\Big ) D_f (x)\,dx\\
&=& \int_\R \int_\R e^{-i\de x\th}\vf(\th) D_f (x)\,d\th dx =\int_\R \vf(\th)\Big (\int_\R e^{-i\de x\th} D_f (x)\,dx\Big )\,d\th\\
&=& \int_{\R}\vf (\th)\widehat{D_f}(\de\th)\,d\th ,
\end{eqnarray*}
thus one gets (\ref{eq3:encad-prob}). Finally combining the latter relation with (\ref{eq1:encad-prob}), (\ref{eq2:encad-prob}), (\ref{eq:four-car}) and (\ref{eq:charX}),
one obtains the lemma. \cqfd

\begin{lemma}
\label{lem:encadexp}
For all real number $u\ge 0$, one has
\begin{equation}
\label{eq1:encadexp}
0\le u-1+e^{-u}\le \frac{u^2}{2}.
\end{equation}
\end{lemma}

\noindent {\sc Proof of Lemma~\ref{lem:encadexp}:} Let $\kappa_1$ and $\kappa_2$ be the functions defined for all real number $u\ge 0$ as, 
$$
\kappa_1(u)=u-1+e^{-u}\mbox{ and } \kappa_2(u)=\frac{u^2}{2}-\kappa_1(u).
$$
One has $\kappa_1(0)=0$ and for all $u\ge 0$, $\kappa_1'(u)=1-e^{-u}\ge 0$; this implies that for each $u\ge 0$, $\kappa_1(u)\ge 0$.
One has $\kappa_2(0)=0$ and for every $u\ge 0$, $\kappa_2 '(u)=\kappa_1(u)\ge 0$; this entails that for each $u\ge 0$, $\kappa_2(u)\ge 0$.
\cqfd

\begin{lemma}
\label{lem:minorh}
Let $h_q$ be the continuous function defined for all $\ga\in [a,b]$, as,
\begin{equation}
\label{eq:defh}
h_q(\ga)=\int_\R |\th|^\ga \vf (\th)\,d\th.
\end{equation}
Then, one has, for all $\ga\in [a,b]$, 
\begin{equation}
\label{eq:encaCphi}
q^{-\ga} h_q(\ga) \le C(\ga) \le q^{\ga} h_q(\ga),
\end{equation}
where $C$ is the continuous and strictly positive function introduced in (\ref{eq3:defC}).
\end{lemma}

\noindent {\sc Proof of Lemma~\ref{lem:minorh}:} First observe that in view of (\ref{eq:sobophi}), the function $h_q$ is well-defined and finite; moreover the dominated convergence 
theorem allows to prove that $h_q$ is continuous on $[a,b]$. Let us now prove that (\ref{eq:encaCphi}) holds. Let $\ga$ be an arbitrary fixed real number belonging 
to the interval $[a,b]$. Assuming that $\a(x)=\gamma$ for all $x\in\R$, then $I(f)$ reduces to a symmetric $\ga$ stable random variable. Next it follows from 
(\ref{eq:expj0}), (\ref{eq1:minoretaj}) and (\ref{eq:LPProbB}), that for all integer $m\ge 1$, one has,
\begin{eqnarray}
\label{eq1:minorh}
\nonumber
&& \int_{\R}\vf (\th)\left(1-\exp\Big\{-\big|q^{-m-1}\th\big|^\ga \int_{\R}\big|f(x)\big|^\ga\,dx\Big\}\right)\,d\th\\
&&\le \P\big (|I(f)|>q^m\big)\\ 
\nonumber
&& \le \int_{\R}\vf (\th)\left(1-\exp\Big\{-\big|q^{-m+1}\th\big|^\ga \int_{\R}\big|f(x)\big|^\ga\,dx\Big\}\right)\,d\th.
\end{eqnarray}
On the other hand, Property 1.2.15 on page 16 in \cite{ST94}, implies that
\begin{equation}
\label{eq2:minorh}
\lim_{m\rightarrow +\infty} q^{m\ga}\, \P\big (|I(f)|>q^m\big) =C(\ga)\int_{\R}\big|f(x)\big|^\ga\,dx.
\end{equation}
Let us now show that 
\begin{equation}
\label{eq3:minorh}
\lim_{m\rightarrow +\infty} q^{m\ga}\int_{\R}\vf (\th)\left(1-\exp\Big\{-\big|q^{-m-1}\th\big|^\ga \int_{\R}\big|f(x)\big|^\ga\,dx\Big\}\right)\,d\th
= q^{-\ga} h_q (\ga)\int_{\R}\big|f(x)\big|^\ga\,dx
\end{equation}
and
\begin{equation}
\label{eq4:minorh}
\lim_{m\rightarrow +\infty} q^{m\ga}\int_{\R}\vf (\th)\left(1-\exp\Big\{-\big|q^{-m+1}\th\big|^\ga \int_{\R}\big|f(x)\big|^\ga\,dx\Big\}\right)\,d\th
= q^\ga h_q (\ga)\int_{\R}\big|f(x)\big|^\ga\,dx.
\end{equation}
We will only prove (\ref{eq4:minorh}) since (\ref{eq3:minorh}) can be obtained similarly. To this end, we will use the dominated convergence theorem. 
It is clear that for all $\th\in\R$,
\begin{equation}
\label{eq5:minorh}
\lim_{m\rightarrow +\infty} q^{m\ga}\vf (\th)\left(1-\exp\Big\{-\big|q^{-m+1}\th\big|^\ga \int_{\R}\big|f(x)\big|^\ga\,dx\Big\}\right)
=q^\ga  |\th|^\ga \vf (\th)\int_{\R}\big|f(x)\big|^\ga\,dx.
\end{equation}
Moreover, it follows from the inequality in the left hand side of (\ref{eq1:encadexp}), that for all integer $m\ge 2$ and real $\th$,
\begin{eqnarray}
\label{eq6:minorh}
\nonumber
q^{m\ga}\big|\vf (\th)\big|\left |1-\exp\Big\{-\big|2^{-m+1}\th\big|^\ga \int_{\R}\big|f(x)\big|^\ga\,dx\Big\}\right|&\le & q^{m\ga}
\big|\vf (\th)\big| \big |q^{-m+1}\th\big|^\ga \int_{\R}\big|f(x)\big|^\ga\,dx\\
&=& q^\ga  |\th|^\ga \big|\vf (\th)\big| \int_{\R}\big|f(x)\big|^\ga\,dx.
\end{eqnarray}
In view of (\ref{eq5:minorh}) and (\ref{eq6:minorh}), we are allowed to apply the dominated convergence theorem, and thus we obtain (\ref{eq4:minorh}).
Finally, putting together (\ref{eq1:minorh}), (\ref{eq2:minorh}), (\ref{eq3:minorh}) and (\ref{eq4:minorh}), one gets (\ref{eq:encaCphi}).
\cqfd

Let us now give some useful properties of the function $T_f$ defined in (\ref{eq1:limdomterm}); easy computations allow to obtain the following two remarks.
\begin{remark}
\label{rem:limdomterm} 
One has for all real number $\xi\ge 1$, 
\begin{equation}
\label{eq2:limdomterm}
\xi^{-b} T_f(1)\le T_f(\xi) \le \xi^{-a} T_f(1).
\end{equation}
\end{remark}

\begin{remark}
\label{rem:ulbTf}
For each real numbers $\de>0$ and $\la>0$,
\begin{itemize}
\item[(i)] when $\de\in (0,1]$, one has
\begin{equation}
\label{eq1:ulbTf}
\de^{-a} T_f (\la)\le T_f(\de\la)\le \de^{-b} T_f(\la),
\end{equation}
\item[(ii)] when $\de >1$, one has 
\begin{equation}
\label{eq2:ulbTf}
\de^{-b} T_f (\la)\le T_f(\de\la)\le \de^{-a} T_f(\la).
\end{equation}
\end{itemize}
\end{remark}

\begin{lemma}
\label{lem:negrest}
For every real number $\xi\ge 1$, we set, 
\begin{equation}
\label{eq1:negrest}
\rho_f(\xi)=\int_\R \big |\vf (\th)\big|\left |\int_\R \big |\xi^{-1} \th f(x)\big|^{\a (x)}\,dx-1+\exp\Big\{-\int_\R \big |\xi^{-1}\th f(x)|^{\a (x)}\,dx\Big\}\right|\,d\th.
\end{equation}
 Then, there is a constant $c(q)>0$, only depending on $b$ and $q$ such that for all $\la \ge q$,
\begin{equation}
\label{eq2:negrest}
\sup_{f\in \S_\a} \frac{\rho_f\big(q^{j_0 (\la,q)+1}\big)}{T_f(\la)}\le c(q)\la^{-a}
\end{equation}
and 
\begin{equation}
\label{eq2bis:negrest}
\sup_{f\in \S_\a} \frac{\rho_f\big(q^{j_0 (\la,q)-1}\big)}{T_f(\la)}\le c(q) \la^{-a};
\end{equation}
recall that $\S_\a$ is the unit sphere of $\F_\a$ and that the strictly positive integer $j_0 (\la,q)$ has been introduced in (\ref{eq:expj0}).
\end{lemma}

\noindent {\sc Proof of Lemma~\ref{lem:negrest}:} It follows from (\ref{eq1:negrest}) and from the right hand side inequality in (\ref{eq1:encadexp}) (in which one takes 
$u=\int_\R \big |\xi^{-1}\th f(x)|^{\a (x)}\,dx)$, that for all real number $\xi\in [1,+\infty)$,
\begin{equation}
\label{eq3:negrest}
\rho_f(\xi)\le \int_\R \Big (\int_\R \big |\xi^{-1} \th f(x)|^{\a (x)}\,dx\Big )^2 \big |\vf (\th)\big|\,d\th.
\end{equation}
Then noticing that 
$$
\Big (\int_\R \big |\xi^{-1}\th f(x)\big|^{\a (x)}\,dx\Big )^2 =\int_{\R}\int_{\R}\big|\xi^{-1}\th f(x_1)\big|^{\a (x_1)}\big|\xi^{-1}\th f(x_2)\big|^{\a (x_2)}\,dx_1 dx_2 ,
$$
and using Fubini-Tonelli Theorem, one gets that
\begin{equation}
\label{eq4:negrest}
\begin{split}
& \int_\R \Big (\int_\R \big |\xi^{-1}\th f(x)|^{\a (x)}\,dx\Big )^2 \big |\vf (\th)\big|\,d\th\\
& =\int_{\R}\int_{\R}\big|\xi^{-1} f(x_2)\big|^{\a (x_1)}\big|\xi^{-1} f(x_2)\big|^{\a (x_2)}
\Big (\int_\R |\th|^{\a (x_1)+\a (x_2)}\big|\vf(\th)\big|\,d\th\Big)\,dx_1 dx_2\\
& \le c_1 \int_{\R}\int_{\R}\big|\xi^{-1} f(x_2)\big|^{\a (x_1)}\big|\xi^{-1} f(x_2)\big|^{\a (x_2)}\,dx_1 dx_2\\
& = c_1 \Big (\int_\R \big |\xi^{-1}f(x)\big |^{\a (x)}\,dx\Big )^2 ,
\end{split}
\end{equation}
where, in view of (\ref{eq:sobophi}), 
$$
c_1=\int_\R \big (1+|\th|\big)^{2b} \big|\vf(\th)\big|\,d\th,
$$ 
is a finite constant. 
Now let us set 
$$
c_2=\max_{y\in [a,b]}\big (C(y)\big)^{-1};
$$
observe that the latter constant is finite since $C$ (see (\ref{eq3:defC})) is a strictly 
positive continuous function on $[a,b]$. Then using (\ref{eq1:limdomterm}), one has that, for all real number $\xi \ge 1$,
\begin{equation}
\label{eq4bis:negrest}
\int_\R \big |\xi^{-1}f(x)\big |^{\a (x)}\,dx\le c_2 T_f (\xi).
\end{equation}
Next, combining (\ref{eq3:negrest}) with (\ref{eq4:negrest}) and (\ref{eq4bis:negrest}), it follows that, 
\begin{equation}
\label{eq5:negrest}
\rho_f (\xi)\le c_3\big (T_f (\xi)\big)^{2}, 
\end{equation}
where $c_3=c_1 c_{2}^2$. Next using (\ref{eq5:negrest}), (\ref{eq:dej0}) and the fact that $T_f$ is a nonincreasing function, one obtains that for all real number $\la \ge q$,
\begin{equation}
\label{eq6:negrest}
0< \frac{\rho_f \big(q^{j_0 (\la,q)+1}\big)}{T_f(\la)}\le c_3 \frac{\big (T_f (q^{j_0 (\la,q)+1})\big)^{2}}{T_f (\la)}\le c_3 T_f (\la)
\end{equation}
and 
\begin{equation}
\label{eq7:negrest}
0< \frac{\rho_f\big(q^{j_0 (\la,q)-1}\big)}{T_f(\la)}\le c_3 \frac{\big (T_f(q^{j_0 (\la,q)-1})\big)^{2}}{T_f(\la)}\le c_3\frac{\big (T_f(q^{-2}\la)\big)^{2}}{T_f(\la)}\le 
q^{4b} c_3 T_f(\la),
\end{equation}
where the latter inequality follows from (\ref{eq1:ulbTf}). Next setting 
$$
c(q)=c_3\Big(1+ q^{4b} \max_{y\in [a,b]} C(y)\Big),
$$ 
and observing that for all $f\in \S_\a$, 
$$
T_f(1)\le \Big(\max_{y\in [a,b]} C(y)\Big)\int_\R \big|f(x)\big|^{\a(x)}\,dx=\Big(\max_{y\in [a,b]} C(y)\Big),
$$ 
then (\ref{eq6:negrest}), (\ref{eq7:negrest})
and Remark~\ref{rem:limdomterm} imply that (\ref{eq2:negrest}) and (\ref{eq2bis:negrest}) hold. 
\cqfd

\begin{lemma}
\label{lem:domterm}
For every real number $\xi\ge 1$, we set,
\begin{equation}
\label{eq1:domterm}
\tau_f(\xi)=\int_\R \vf(\th)\Big (\int_\R \big |\xi^{-1}\th f(x)\big|^{\a (x)}\,dx\Big )\,d\th.
\end{equation}
Then, one has, for all real number $\xi \ge q$,
\begin{equation}
\label{eq3:domterm}
T_f( q \xi)\le \tau_f(\xi) \le T_f\big( q^{-1} \xi\big);
\end{equation}
moreover
\begin{equation}
\label{eq4:domterm}
q^{-2b}\le \inf_{\la\in [q,+\infty)} \frac{\tau_f\big(q^{j_0(\la,q)+1}\big)}{T_f(\la)}\le \sup_{\la\in [q,+\infty)} \frac{\tau_f \big(q^{j_0(\la,q)-1}\big)}{T_f(\la)}\le q^{3b}.
\end{equation}
\end{lemma}

\noindent {\sc Proof of Lemma~\ref{lem:domterm}:} First observe that (\ref{eq:sobophi}) and the fact that $f\in\F_\a$, imply that for all real number $\xi\ge 1$,
$$
\int_\R \int_\R \big|\vf(\th)\big|\big |\xi^{-1}\th f(x)\big|^{\a (x)}\,dx\Big )\,dx d\th\le \Big (\int_\R \big (1+|\th|\big)^{b} \big|\vf(\th)\big|\,d\th\Big)
\times \Big ( \int_{\R} \big|f(x)\big|^{\a(x)}\,dx\Big )<\infty.
$$
Therefore, we are allowed to use Fubini Theorem and we obtain that, for all real number $\xi\ge 1$,
\begin{equation}
\label{eq5:domterm}
\tau_f (\xi) = \int_\R \big |\xi^{-1}f(x)\big |^{\a (x)} h_q \big(\a(x)\big)\,dx.
\end{equation}
where the function $h_q $ has been introduced in~(\ref{eq:defh}). Next it follows from (\ref{eq:encaCphi}) that, for all real number $\xi\ge 1$,
$$
\int_\R \big |\xi^{-1}f(x)\big |^{\a (x)} q^{-\a(x)} C\big(\a(x)\big)\,dx\le \tau_f(\xi)\le \int_\R \big |\xi^{-1}f(x)\big |^{\a (x)} q^{\a(x)} C\big(\a(x)\big)\,dx;
$$
thus, in view of (\ref{eq1:limdomterm}), one gets (\ref{eq3:domterm}). Let us now prove that the first inequality in (\ref{eq4:domterm}) holds. Assume that the real $\la\ge q$
is arbitrary. It follows from the first inequality in (\ref{eq3:domterm}), (\ref{eq:dej0}), the fact that $T_f$ is a nonincreasing function and the first inequality in (\ref{eq2:ulbTf}) (in which one takes $\de=q^2$), that
$$
\frac{\tau_f \big(q^{j_0(\la,q)+1}\big)}{T_f(\la)}\ge \frac{T_f\big(q^{j_0(\la,q)+2}\big)}{T_f(\la)}\ge \frac{T_f\big(q^{2}\la \big)}{T_f(\la)}\ge q^{-2b};
$$
thus we obtain the first inequality in (\ref{eq4:domterm}). Next observe that (\ref{eq5:domterm}) and the fact that $h_q$ is a (strictly) positive function 
(this is a straightforward consequence of (\ref{eq:encaCphi})) imply that $\tau_f$ is a nonincreasing function, which in turn entails 
that for all real number $\la\ge q$ one has $\tau_f\big(q^{j_0(\la,q)+1}\big)\le\tau_f \big(q^{j_0(\la,q)-1}\big)$ and, as consequence, that the second inequality in (\ref{eq4:domterm}) is satisfied. Finally, let us show that the last inequality in (\ref{eq4:domterm}) holds. It follows from the second inequality in (\ref{eq3:domterm}), (\ref{eq:dej0}), the fact that $T_f$ is a nonincreasing function and the second inequality in (\ref{eq1:ulbTf}) (in which one takes $\de=q^{-3}$), that
$$
\frac{\tau_f \big(q^{j_0(\la,q)-1}\big)}{T_f(\la)}\le\frac{T_f \big(q^{j_0(\la,q)-2}\big)}{T_f(\la)}\le \frac{T_f\big(q^{-3}\la\big)}{T_f(\la)} \le q^{3b};
$$
thus we obtain the last inequality in (\ref{eq4:domterm}).
\cqfd

\begin{lemma}
\label{lem:minoretaj}
Let $\eta_f$ be the function introduced in (\ref{eq1:minoretaj}) and let $c(q)>0$ be the constant introduced in Lemma~\ref{lem:negrest}. Then one has
for all $f\in \S_\a$,  and for all real number $\la\ge q$ 
\begin{equation}
\label{eq2:minoretaj}
q^{-2b}-c(q)\la^{-a}\le \frac{\eta_f\big(q^{j_0(\la,q)+1}\big)}{T_f(\la)}\le \frac{\eta_f \big(q^{j_0(\la,q)-1}\big)}{T_f(\la)}\le q^{3b}+c(q)\la^{-a}.
\end{equation}
\end{lemma}

\noindent {\sc Proof of Lemma~\ref{lem:minoretaj}:} In view of (\ref{eq1:minoretaj}), for 
all real number $\xi\ge 1$, one can write,
\begin{eqnarray*}
\eta_f (\xi) &=& \int_\R \vf(\th)\Big (\int_\R \big |\xi^{-1}\th f(x)\big|^{\a (x)}\,dx\Big )\,d\th\\
&&+\int_\R \vf (\th)\left (1-\int_\R \big |\xi^{-1} \th f(x)\big|^{\a (x)}\,dx-\exp\Big\{-\int_\R \big |\xi^{-1}\th f(x)|^{\a (x)}\,dx\Big\}\right)\,d\th;
\end{eqnarray*}
thus, it follows from the triangle inequality (\ref{eq1:domterm}) and (\ref{eq1:negrest}), that for 
each real number $\xi\ge 1$, one has,
\begin{equation}
\label{eq3:minoretaj}
\tau_f(\xi)-\rho_f(\xi)\le \eta_f(\xi) \le \tau_f(\xi)+\rho_f(\xi).
\end{equation}
Let us now prove that the first inequality in (\ref{eq2:minoretaj}) holds. Using the first inequality in (\ref{eq3:minoretaj}) as well as the first inequality in (\ref{eq4:domterm}), one has for all real number $\la\ge q$,
$$
q^{-2b}-\frac{\rho_f\big(q^{j_0(\la,q)+1}\big)}{T_f(\la)}\le \frac{\eta_f\big(q^{j_0(\la,q)+1}\big)}{T_f(\la)};
$$
then the first inequality in (\ref{eq2:minoretaj}) results from (\ref{eq2:negrest}). Next, observe that (\ref{eq:LPProbB}) clearly implies that the second inequality in (\ref{eq2:minoretaj}) is satisfied. Finally, let us show that the last inequality in (\ref{eq2:minoretaj}) holds. Using the second inequality in (\ref{eq3:minoretaj}) as well as the second inequality in (\ref{eq4:domterm}), one has for all real number $\la\ge q$,
$$
\frac{\eta_f\big(q^{j_0(\la,q)-1}\big)}{T_f(\la)}\le q^{3b}+\frac{\rho_f\big(q^{j_0(\la,q)-1}\big)}{T_f(\la)};
$$
then the last inequality in (\ref{eq2:minoretaj}) results from (\ref{eq2bis:negrest}).
\cqfd
\\

Now we are in position to prove Theorem~\ref{th:main}.

\noindent {\sc Proof of Theorem~\ref{th:main}:} Let $c(q)$ be the constant introduced in Lemma~\ref{lem:negrest}. Using (\ref{eq:LPProbB}) and (\ref{eq2:minoretaj}), one has for all $f\in\S_\a$ and all real number $\la\ge q$,
$$
 q^{-2b}-c(q)\lambda^{-a}\le \frac{\eta_f\big(q^{j_0(\la,q)+1}\big)}{T_f(\la)}\le \frac{\P\big(|I(f)|>\la\big)}{T_f(\la)}\le  \frac{\eta_f\big(q^{j_0(\la,q)-1}\big)}{T_f(\la)}\le 
 c (q) \lambda^{-a}+q^{3b},
$$
which implies that,
$$
q^{-2b}-1-c(q)\lambda^{-a}\le \frac{\P\big(|I(f)|>\la\big)}{T_f(\la)}-1\le c(q)\lambda^{-a}+q^{3b}-1
$$
and, as a consequence, that
$$
\sup_{f\in\S_\a}\Big|\frac{\P\big(|I(f)|>\la\big)}{T_f(\la)}-1\Big|\le |q^{-2b}-1|+|q^{3b}-1|+c(q)\lambda^{-a}.
$$
One has therefore, 
$$
\limsup_{\la\rightarrow +\infty} \left\{\sup_{f\in\S_\a}\Big|\frac{\P\big(|I(f)|>\la\big)}{T_f(\la)}-1\Big|\right\}\le |q^{-2b}-1|+|q^{3b}-1|.
$$
Finally letting $q>1$ goes to $1$, one obtains the theorem. \cqfd

\bibliographystyle{plain}
\begin{small}

\end{small}

\end{document}